\input amstex
\documentstyle{amsppt}

\magnification=\magstep1

\define\a{\alpha}

\define\la{\lambda}

\define\conv{\text{\rm conv}}
\define\vo #1{\text{\rm vol} \left ( #1 \right )}
\define\etc{, \dots ,}

\define\rn{$\Bbb R^n \,$}
\define\RM{$\Bbb R^m \,$}
\define\Rm{\Bbb R^m}
\define\Rn{\Bbb R^n}

\define\nor #1{\left \| #1 \right \|}
\define\enor #1{\Bbb E \, \nor{#1}}
\define\tens #1{#1 \otimes #1}
\define\pr#1#2{\langle {#1} , {#2} \rangle}

\topmatter

\title
Sections of the difference body
\endtitle

\author
M. Rudelson
\endauthor

\thanks
This research was  supported in part  by NSF grant DMS-9706835.
\endthanks

\affil
Texas A \& M University
\endaffil

\address
Department of Mathematics, \hfill \break
Texas A \& M University, \hfill \break
College Station, TX 77843.
\endaddress

\email
Mark.Rudelson\@math.tamu.edu
\endemail

\vskip 1in

\abstract
Let $K$ be an $n$-dimensional convex body. 
Define the difference body by 
$$
K-K= \{ x-y \mid x,y \in K \}.
$$
We estimate the volume of the section of $K-K$ by a linear subspace $F$ via
the maximal volume of sections of $K$ parallel to $F$. 
We prove that for any $m$-dimensional subspace $F$ there exists $x \in \Rn$,
such that 
$$
\vo{(K-K) \cap F} \le 
C^m  \left ( \min \left ( \frac{n}{m}, \sqrt{m} \right ) \right )^m
\cdot  \vo{K \cap (F+x)},
$$
for some absolute constant $C$.
We show that for small dimensions of $F$ this estimate is exact up to 
a multiplicative constant.
\endabstract

\endtopmatter

\vskip .2in
\document

\head 1. Introduction. \endhead

Let $K$ be an $n$-dimensional convex body.
Define the difference body by 
$$
K-K= \{ x-y \mid x,y \in K \}.
$$
In 1957 Rogers and Shephard \cite{R-S} proved that
$$
\vo{K-K} \le {{2n} \choose n} \vo{K}.
$$
A simpler proof was found later by Chakerian \cite{C}.

Let $F$ be an $m$-dimensional linear subspace of \rn and let $P_F$ be the
orthogonal projection onto $F$.
It follows from the inequality of Rogers and Shephard that 
$$
\left ( \frac{\vo{P_F(K-K)}}{\vo{P_F K}} \right )^{1/m} 
\le {{2m} \choose m}^{1/m} < 4 .
$$
Here and later we  denote by vol the volume in the relevant dimension.

For some problems it would be interesting to obtain a similar
estimate for the volumes of {\it sections} of $K-K$.
In particular, would the  expression
$$
R(K,F) = \left ( \frac{\vo{(K-K) \cap F}}{\sup_{x \in \Rn} \vo{K \cap (F+x)}} 
\right )^{1/m}
$$
be uniformly bounded?

Although, as it is shown below, the answer to this question is negative,
some estimates of this ratio are possible.
Our main result is the following
\proclaim{Theorem 1}
Let $K \subset \Rn$ be a convex body and let $F \subset \Rn$ be 
an $m$-dimensional subspace. 
Then
$$
\vo{(K-K) \cap F} \le C^m  \varphi^m (m,n) \cdot  \sup_{x \in \Rn}
\vo{K \cap (F+x)},
$$
where
$$
\varphi (m,n) = \min \left ( \frac{n}{m}, \sqrt{m} \right ).
$$
\endproclaim
Here and later $C$ denotes an absolute constant whose value may change 
from line to line.

This result can be applied to estimating the Banach -- Mazur distance between 
two non-symmetric convex bodies. 
To use random rotations for such an estimate one has to put the bodies into
some specific positions. 
This can be achieved by comparison of the positions of the difference body 
and the body itself.
We are going to present the details in a separate paper.

 It follows from Theorem 1 that $R(K,F)$ is bounded for  $m$ proportional
to $n$ and for a small $m$.
This suggests that  $R(K,F)$ should be bounded for all dimensions.
Surprisingly, this is not the case.
Namely, the following Theorem implies that for some body $K \subset \Rn$
and $F \subset \Rn, \text{dim}(F)=m$
$$
R(K,F) \ge c \sqrt{\log n}
$$ 
when $c \log n \le m \le n^{\a}$ and $\a \in (0,1)$.

\proclaim{Theorem 2}
For any $m<n$ there exists a convex body $K \subset \Rn$ and a subspace
$F \subset \Rn$ of dimension $m$ such that for any $x \in \Rn$
$$
\vo{(K-K) \cap F} \ge C^m \psi^m (m,n) \cdot \vo{K \cap (F+x)},
$$
where 
$$
\psi(m,n)= \min \left (
\sqrt{\log \left ( \frac{n}{m}+1 \right )}, \sqrt{m} \right ).            
$$
\endproclaim

 Notice that Theorem 2 implies that the estimate obtained in Theorem 1 
is exact for $m \le c \log n$.

\subhead Acknowledgment \endsubhead
The author thanks F.~Barthe and the referees for their valuable remarks.

\head 2.  Upper estimate. \endhead

The proof of Theorem 1 consists of two steps. 
First we reduce the problem to a question of comparing the volume of 
projection and the volume of parallel sections of a certain convex body.
Then we use the Rogers -- Shephard inequality and the John
decomposition to complete the proof.

Denote by $V_m(D)$ the $m$-th intrinsic volume of a body $D$ \cite{S}.
Consider the following integral
$$
I(K,F) = \int_F V_m \big ( K \cap (K+x) \big ) dx.
$$
To prove the Theorem we shall estimate $I(K,F)$ from above and from below.

For the lower estimate
we apply the following Lemma due to Chakerian \cite{C}.
\proclaim{Lemma 1}
Let $B \subset \Rm, \ 0 \in B$ be a convex body. 
Let $h: B \to \Bbb R$ be a
non-negative concave function and let $f: \Bbb R \to \Bbb R$ be increasing.
Then
$$
\int_B f(h(x)) dx \ge 
m \cdot \vo{B} \cdot \int \limits_0^1 f(t \cdot h(0)) (1-t)^{m-1} dt. \qquad 
\qed
$$
\endproclaim

For $x \in F$ let
$$
h(x)=V_m^{1/m}  \big (  D_x(K) \big ), \qquad f(t)=t^m,
$$
where
$$
D_x(K)= K \cap (K+x)
$$
It follows from the Alexandrov -- Fenchel inequality that the 
intrinsic volumes satisfy the General Brunn -- Minkowski inequality.
Namely,  for any two 
bodies $B,D$ and for any number $0 \le \la \le 1$
$$
V_m^{1/m} ( \la B + (1-\la) D ) \ge  \la V_m^{1/m}(B) + (1-\la ) V_m^{1/m}(D)
\tag 1
$$
\cite{S, Th. 6.4.3, p.339}.
Since for any $x, \bar x$
$$
\la D_x(K) + (1-\la) D_{\bar x} (K) \subset D_{ \la x + (1-\la) \bar x} (K),
$$
it follows from (1) that $h(x)$ is a concave function.
By  Lemma 1,
$$
\align
\int_{(K-K) \cap F} h^m(x) \, dx 
&\ge m \cdot \vo{(K-K) \cap F} \cdot 
\int_0^1 \big ( t \cdot h(0) \big )^m \cdot (1-t)^{m-1} \, dt \\
&= \vo{(K-K) \cap F} \cdot h^m(0) \cdot {{2m} \choose m}^{-1}.
\endalign
$$
So, we get that 
$$
\aligned
\vo{(K-K) \cap F} &\le {{2m} \choose m} \cdot V_m^{-1} (K) \cdot 
\int_{(K-K) \cap F} V_m \big ( K \cap (K+x) \big ) \, dx \\   
&\le 4^m \cdot V_m^{-1} (K) \cdot 
\int_F V_m \big ( K \cap (K+x) \big ) \, dx.
\endaligned  \tag 2 
$$

To estimate $I(K,F)$ we apply Crofton's formula 
\cite{S, formula (4.5.9), p. 235}.
 Let $\Bbb A (n,n-m)$ be the set of all 
$(n-m)$-dimensional affine subspaces of \rn and let $\mu$ be the Haar 
measure on $\Bbb A (n,n-m)$.
By  Crofton's formula, we get 
$$
V_m ( K \cap (K+x)) =
C_{n,m} \cdot  \int_{\Bbb A (n,n-m)}
\chi (K \cap (K+x) \cap E) \, d \mu (E),
$$
where $C_{n,m}$ is a constant depending on $n$ and $m$.
By Fubini's theorem, 
$$
\align
I(K,F) &=\int_{(K-K) \cap F}  V_m ( K \cap (K+x)) \, dx \\
&= C_{n,m} \cdot \int_F \int_{\Bbb A (n,n-m)}
\chi (K \cap (K+x) \cap E) \, d \mu (E) \, dx \\
&= C_{n,m} \cdot \int_{\Bbb A (n,n-m)} \text{mes } \{ x \in F \, | \,
(K+x) \cap (K \cap E) \neq \emptyset \}  \, d \mu (E),  \tag 3
\endalign  
$$
where mes is the Lebesgue measure on $F$.
Let $\Bbb A_F$ be the set of all $(n-m)$-dimen\-sional affine subspaces 
which are transversal to $F$:
$$
\Bbb A_F = \{ E \in \Bbb A (n,n-m) \ | \ \text{card }(E \cap F) = 1 \}.
$$
Since $\mu ( \Bbb A (n,n-m) \backslash \Bbb A_F) =0$, we can integrate in (3) 
only over $\Bbb A_F$.
Then (3) can be estimated above by 
$$
\align
&C_{n,m} \cdot \int_{\Bbb A_F} \chi (K \cap E) \, d \mu (E) \cdot
\sup_{E \in \Bbb A_F} \text{mes } \{ x \in F \, | \,
(K+x) \cap (K \cap E) \neq \emptyset \} \\
&= V_m (K) \cdot \sup_{E \in \Bbb A_F} \text{mes } \{ x \in F \, | \,
(K+x) \cap (K \cap E) \neq \emptyset \}.
\endalign
$$

To complete the proof of Theorem 1 we have to prove the following 
\proclaim{Claim}
For any $m$-dimensional linear subspace $F \subset \Rn$ and
any $(n-m)$-dimen\-sional affine subspace $E \subset \Rn$, such that
$E$ and $F$ intersect at one point only,
$$
\text{\rm mes } \{ x \in F \, | \, (K+x) \cap (K \cap E) \neq \emptyset \}
\le C^m \varphi^m (m,n) \cdot \sup_{y \in \Rn} \vo{(K+y) \cap F}.
$$
\endproclaim
\demo{Proof of the Claim}
Since the statement of the Claim is invariant under translations, we may 
assume that $E \cap F = \{ 0 \}$.
Also, let $T: \Rn \to \Rn$ be an invertible linear operator, such that 
$T|_F = id$ and $T|_E = F^{\perp}$. The Claim is invariant under $T$, 
so we may assume that $E$ and $F$ are orthogonal.

Define 
$$
Z=K \cap \big ( (K \cap E ) +F \big ).
$$
Let $P_E,P_F$ be orthogonal projections onto $E$ and $F$ respectively.
We have
$$
\align
&\text{mes } \{ x \in F \, | \, (K+x) \cap (K \cap E) \neq \emptyset \} \\
&=\text{mes } \{ x \in F \, | \, 
K \cap \big ( (K \cap E)-x \big ) \neq \emptyset \} \\
&=\text{mes }\Big ( P_F \big ( (K \cap (K \cap E)-F) \big ) \Big )
=\vo{P_F (Z)}.
\endalign
$$
By the construction of $Z$ we have
$$
Z \cap E \subset P_E Z \subset P_E \big ( (K \cap E ) +F \big )
=K \cap E = Z \cap E. \tag 4
$$
Since $Z \subset K$, and ${n \choose m} \le e^m (n/m)^m$, it is enough
 to prove that
$$
\vo{P_F (Z)} \le  {n \choose m} \cdot \sup_{y \in E} \vo{(Z+y) \cap F} \tag i
$$
and
$$
\vo{P_F (Z)} \le C^m  m^{m/2} \cdot \sup_{y \in E} \vo{(Z+y) \cap F}. \tag ii
$$
\enddemo
\demo{Proof of (i)}
By (4),
$$
\aligned
\vo{Z} &\le \vo{P_E Z} \cdot \sup_{y \in E} \vo{(Z+y) \cap F} \\
&=\vo{Z \cap E} \cdot \sup_{y \in E} \vo{(Z+y) \cap F}.
\endaligned \tag 5
$$
From the other side, another inequality of Rogers and Shephard \cite{R-S}
implies that
$$
\vo{Z} \ge {n \choose m}^{-1} \cdot \vo{P_F Z} \cdot \vo{Z \cap E}. \tag 6
$$
Now (i) follows from the combination of (5) and (6).
\qed
\enddemo

\demo{Remark}
Using the inequality (6) of Rogers and Shephard in the proof of (i) leads to a 
gap between the upper and lower estimates of $\varphi (m,n)$.
Although the  Rogers and Shephard inequality is exact, it holds as an equality
for the bodies of the form $Z=\conv (Z \cap E, Z \cap F)$, while for such 
bodies $P_F (Z) =Z \cap F$.
\enddemo

\demo{Proof of (ii)}
Without loss of generality we may assume that the ellipsoid of minimal volume
containing $P_F Z$ is $B_2^m$. 
Then there exists a John's decomposition of the identity operator.
Namely, there exist 
$M \le (n+3)n/2$  contact points $x_1 ,\dots, x_M \in S^{m-1} \cap Z$ and 
$M$ positive numbers $c_1, \dots, c_M$ satisfying the following 
system of equations
$$
\align
id &= \sum_{i=1}^M c_i \, x_i \otimes x_i \\
0 &= \sum_{i=1}^M c_i \, x_i. 
\endalign
$$
Here by $id$ we denote the identity operator in \RM.

Since $x_i \in P_F Z$, we can  choose the points $y_i \in P_E Z$ so that 
$x_i+y_i \in Z$.
Define 
$$
u=  \sum_{i=1}^M \frac{c_i}{m} \, y_i.
$$
Since
$$
\sum_{i=1}^M \frac{c_i}{m} =1,
$$
$u \in P_E Z = Z \cap E$.
Notice that $y_1 \etc y_M \in Z$.
So,
$$
\align
Z \cap (F+u) &\supset 
\sum_{i=1}^M \frac{c_i}{m} \cdot \big ( Z \cap (F+y_i) \big ) \\
&\supset \sum_{i=1}^M \frac{c_i}{m} \cdot [ y_i, y_i + x_i ] 
= \left ( \sum_{i=1}^M \frac{c_i}{m} \cdot [ 0, x_i ] \right ) + u \\
&= \frac {1}{2} \cdot
\left ( \sum_{i=1}^M \frac{c_i}{m} \cdot [ -x_i, x_i ] \right ) + u.
\endalign
$$
Here $\sum$ means the Minkowski sum and $[x,y]$ denotes the segment
joining $x$ and $y$.
Put 
$$
W=  \sum_{i=1}^M \frac{c_i}{m} \cdot [ -x_i, x_i ].
$$
Then, by \cite{B, Lemma 4}, we have
$$
\align
\vo{W} &\ge 2^m m^{-m},
\intertext{so}
\vo{P_F Z} &\le \vo{B_2^m} \le C^m m^{m/2} \cdot \vo{W}. \qed
\endalign
$$

\enddemo

Notice that $\varphi(m,n) \le n^{1/3}$.
So, we have the following immediate
\proclaim{Corollary}
Let $K \subset \Rn$ be a convex body and let $F \subset \Rn$ be 
an $m$-dimensional subspace. 
Then
$$
\vo{(K-K) \cap F} \le \left (C \cdot  n^{1/3} \right )^m \cdot  
\sup_{x \in \Rn} \vo{K \cap (F+x)}.  \qed
$$
\endproclaim

\head 3. Lower estimate. \endhead

We now turn to the proof of Theorem 2.

Assume first that  $ n-m+1 \ge 5^m$.
In this case we have to prove Theorem 2 for $\psi(m,n)= \sqrt{m}$.
The assumption guarantees that one can find  points $z_1 \etc z_{n-m+1}$ 
on the unit sphere of $F$ which form a $(1/2)$-net.
Let 
$j_1 \etc j_{n-m+1}$ be the vertices of the standard simplex in the 
space $F^{\perp}$.        
Put
$$
K=\conv (j_1 \pm z_1 \etc j_{n-m+1} \pm z_{n-m+1}).
$$
Since 
$$
(K-K) \cap F \supset 2 \conv (\pm z_1 \etc \pm z_{n-m+1}) \supset B_2^m,
$$ 
we have to prove that for any $x \in \conv (j_1 \etc j_{n-m+1})$
$$
\vo{K \cap (F+x)} \le \left ( \frac{c}{\sqrt{m}} \right)^m \vo{B_2^m}.
$$
Assume that 
$$
\align
x &= \sum_{i=1}^{n-m+1} \lambda_i j_i,
\intertext{where}
\lambda_i &\ge 0 \text{ and } \sum_{i=1}^{n-m+1} \lambda_i =1.
\endalign
$$
Then 
$$
K \cap (F+x) = \sum_{i=1}^{n-m+1} \lambda_i [j_i-z_i, j_i+z_i]=
\sum_{i=1}^{n-m+1} \lambda_i [-z_i, z_i] +x.
$$
Let $T_1 \etc T_N \in O(m)$ be random rotations in \RM. 
By the Brunn -- Minkowski inequality 
$$
\align
\vo{K \cap (F+x)} 
&\le \vo{\frac{1}{N} \sum_{s=1}^{N} T_s ((K-x) \cap F)} \\
&=\vo{\frac{1}{N}  \sum_{s=1}^{N} 
\sum_{i=1}^{n-m+1} \lambda_i [-T_s z_i, T_s z_i]}. \tag 7
\endalign
$$
For a sufficiently large $N$ 
$$
\frac{1}{N}  \sum_{s=1}^{N} \lambda_i [-T_s z_i, T_s z_i]
\subset \frac{2}{\sqrt{m}} B_2^m
$$
so (7) does not exceed
$$
\vo{\frac{2}{\sqrt{m}} \sum_{i=1}^{n-m+1} \lambda_i B_2^m}
=\vo{\frac{2}{\sqrt{m}}  B_2^m}.
$$

Now assume that $n-m+1<5^m$ and let $k$ be the largest integer such that
$5^k \cdot (m/k) \le n-m+1$.
Since in this case $k \le c \log ( n/m + 1)$, it is enough to prove Theorem~2
for $\psi(m,n)=\sqrt{k}$. 
We shall use a construction which is similar to \cite{F-J, p. 96--97}.
Assume for simplicity that $L=m/k$ is an integer.
Let $e_1 \etc e_m$ be an orthonormal basis of $F$.
For $l=1 \etc L$ put
$$
F_l= \text{span } \{ e_i \mid i=k(l-1)+1 \etc kl \}.
$$
Let $z_1^l \etc z_{N_l}^l$ be an $1/2$-net on the unit sphere of $F_l$.
Since $5^k \cdot (m/k) \le n-m+1$, we may assume that the total number 
of elements in these nets is $n-m+1$.
Let us reorder the sequences $\{z_i^l \}$ into one sequence 
$\{z_i \}_{i=1}^{n-m+1}$.
Let $j_1 \etc j_{n-m+1}$ be the vertices of the standard simplex 
in $F^{\perp}$.
Define as before
$$
K=\conv (j_1 \pm z_1 \etc j_{n-m+1} \pm z_{n-m+1}).
$$
Then we have 
$$
(K-K) \cap F \supset 2 \conv (\pm z_1 \etc \pm z_{n-m+1}).
$$ 
Since the sequence $z_1 \etc z_{n-m+1}$ contains the $(1/2)$-nets for 
the unit spheres of the spaces $F_l$,
$$
(K-K) \cap F \supset \conv (B_2^m \cap F_1 \etc B_2^m \cap F_L).
$$
Put $B_l=B_2^m \cap F_l$.
We have to prove that for any $x \in \conv (j_1 \etc j_{n-m+1})$
$$
\vo{K \cap (F+x)} \le \left ( \frac{c}{\sqrt{k}} \right )^m \cdot
\vo{\conv (B_1 \etc B_L)}. \tag 8
$$
Assume that 
$$
\align
x &= \sum_{i=1}^{n-m+1} \lambda_i j_i,
\intertext{where}
\lambda_i &\ge 0 \text{ and } \sum_{i=1}^{n-m+1} \lambda_i =1.
\endalign
$$
Then as before we have 
$$
K \cap (F+x) = 
\sum_{i=1}^{n-m+1} \lambda_i [-z_i, z_i] +x.
$$
Let $T_1^l \etc T_M^l: F_l \to F_l$, be random  rotations of  $F_l$.
Denote by $I_l$  the set of indexes $i$ for which $z_i \in F_l$.
Then
$$
\frac{1}{M} \sum_{s=1}^M \sum_{i \in I_l} 
\lambda_i [-T_s^l z_i, T_s^l z_i] \subset \frac{2}{\sqrt{k}} \mu_l B_l.
$$
Here $\mu_l = \sum_{i \in I_l} \lambda_i$, so
$$
\sum_{l=1}^L \mu_l =1.
$$
Arguing as before, we prove that
$$
\align
\vo{K \cap (F+x)} 
&\le \vo{ \sum_{l=1}^L \frac{1}{M} \sum_{s=1}^M  
\sum_{i \in I_l} \lambda_i [-T_s^l z_i, T_s^l z_i]} \\
&\le \vo{ \sum_{l=1}^L \frac{2}{\sqrt{k}} \mu_l  B_l}\\
&=\left ( \frac{2}{\sqrt{k}} \right )^m \cdot  
\left (  \prod_{l=1}^L \mu_l \right )^k \cdot \left ( \vo{B_1} \right )^L.   
\endalign 
$$
By the inequality between the arithmetic and the geometric mean,
$$
\vo{K \cap (F+x)} \le \left ( \frac{2}{\sqrt{k}} \right )^m \cdot L^{-kL}
\cdot \left ( \vo{B_1} \right )^L.   \tag 9
$$
To complete the proof we apply the following easy Lemma, 
which can be proved by induction.
\proclaim{Lemma 3}
Let $\Bbb R^{kL}=F_1 \oplus \ldots \oplus F_L$, where $F_1 \etc F_L$ are 
mutually orthogonal subspaces of dimension $k$. Let $B_l=B_2^{kL} \cap F_l$.
Then 
$$
\vo{\conv (B_1 \etc B_L)} = \frac{(k!)^L}{(kL)!} \Big [ \vo{B_1} \Big ]^L. 
\qquad \qed
$$
\endproclaim
\demo{Remark}
A generalization of this formula appears in \cite{M}.
\enddemo
Since 
$$
\frac{(k!)^L}{(kL)!} \le C^L \cdot L^{kL} = C^L \cdot L^m
$$
the inequality (8) follows from (9) and the Lemma.

\vskip .5in
\head References \endhead

\item{[B]}
K. Ball,
Shadows of convex bodies,
{\it Trans of AMS} {\bf 327}, no. 2,  (1991),
891--901.

\item{[C]}
G. Chakerian,
Inequalities for the difference body of a convex body,
{\it Proc. of AMS}  {\bf  18}  (1967),
 879--884.

\item{[F-J]}
 T. Figiel and W.B. Johnson,
Large subspaces of $\ell_{\infty}^n$ and estimates of the Gordon --
Lewis constants,
{\it Israel J. Math.}  {\bf 37} (1980),
92--112.

\item{[M]}
P. McMullen, 
The volume of certain convex sets,
{\it Math. Proc. Camb. Phil. Soc.} {\bf 1991} (1982),
91.

\item{[R-S]}
C.A. Rogers and G.C. Shephard,
The difference body of a convex body
{\it Arch. Math.} {\bf 8} (1957),
220--233.

\item{[S]}
 R. Schneider,
Convex bodies: the Brunn -- Minkowski theory,
 Cambridge University Press,
 Cambridge,  1993.

\enddocument